\documentstyle [12pt]{article}

\title {On uncountable hypersimple unidimensional theories}
\author {Ziv Shami\\Ariel University}
\newtheorem {theorem}{Theorem}[section]

\newtheorem {definition}[theorem]{Definition}

\newtheorem {fact}[theorem]{Fact}

\newtheorem {corollary}[theorem]{Corollary}

\newtheorem {proposition}[theorem]{Proposition}

\newtheorem {subclaim}[theorem]{Subclaim}

\def\proof {\noindent \textbf{Proof:} }

\newsavebox{\indbin}
\savebox{\indbin}{\begin{picture}(0,0)
\newlength{\gnu}
\settowidth{\gnu}{$\smile$} \setlength{\unitlength}{.5\gnu} \put(-1,-.65){$\smile$}
\put(-.25,.1){$|$}
\end{picture}}
\newcommand{\nonfork}[3]
{\mbox{$\begin{array}{ccc} \mbox{$#1$} & \usebox{\indbin} & \mbox{$#2$} \\
        & \mbox{$#3$} &
\end{array}$}}
\newcommand{\nonforkempty}[2]
{\mbox{$\begin{array}{ccc} \mbox{$#1$} & \usebox{\indbin} & \mbox{$#2$}
\end{array}$}}

\newcommand{\forkempty}[2]
{\mbox{$\begin{array}{ccc} \mbox{$#1$} & \!\mbox{$\!\!\not\!\:\usebox{\indbin}$} & \mbox{$#2$}
\end{array}$}}

\def\card #1 {{\vert #1 \vert}}

\def\CC {{\cal C}}
\def\RR {{\cal R}}

\def\UU {{\cal U}}

\def\WW {{\cal W}}

\begin{document}
\maketitle

\begin{abstract}
We extend the dichotomy between 1-basedness and supersimplicity proved in [S1]. The generalization
we get is to arbitrary language, with no restrictions on the topology (we do not demand
type-definabilty of the open set in the definition of essential 1-basedness from [S1]). We conclude
that every (possibly uncountable) hypersimple unidimensional theory that is not s-essentially
1-based by means of the forking topology is supersimple. We also obtain a strong version of the
above dichotomy in the case where the language is countable.
\end{abstract}

\section{Introduction}
Shelah has defined unidimensional theories as (stable) theories in which any two sufficiently large
$\vert T\vert^+$-saturated models of the same cardinality are isomorphic. For stable theories this
definition is equivalent to the requirement that any two non-algebraic types are non-orthogonal.
This requirement serves as the definition of unidimensionality for the larger class of simple
theories. A problem posed by Shelah was whether any unidimensional stable theory is superstable.
Around 1986 Hrushovski has solved the problem by answering it in the affirmative [H1].

Several years after the discovery of Kim [K] that the algebraic properties of forking (symmetry and
transitivity) can be proved for simple theories (1996) and the development of the basic machinery
[K,KP,HKP], there were several attempts to generalize the above result of Hrushovski to the simple
case. A generalization of this proof along the the same lines seems very problematic because of the
lack of definability of types, and so many of the results on definable groups in stable theories do
not seem to generalize to simple theories in a direct way.

In 2003, we observed that any small simple unidimensional theory is supersimple [S3]. A bit later,
Pillay [P] has proved that any countable hypersimple theory (i.e. a simple theory that eliminates
hyperimaginaries) with the wnfcp (the weak non finite cover property) is supersimple; this proof
builds on ideas from Hrushovski's old proof of the result for countable stable theories [H0] and
some machinery from the theory of lovely-pairs [BPV]. This has been extended by Pillay [P1] to any
countable low hypersimple theory using the result on elimination of the "there exists infinitely
many" quantifier [S2]. In 2008, it has been proved that any countable hypersimple unidimensional
theory is supersimple [S1]. An important notion that used in [S1] is the forking topology (or the
$\tau^f$-topology); this is a variant of the topology used in [H0] and [P]: for variables $x$ and
set $A$ the forking topology on $S_x(A)$ is defined as the topology whose basis is the collection
of all sets of the form $\UU=\{a \vert \phi(a,y)$ forks over $A \}$, where $\phi(x,y)\in L(A)$.

The goal of this paper is to reduce the problem on supersimplicity of general hypersimple
unidimensional theories (possibly uncountable) to the case where the theory is s-essentially
1-based by means of the forking topology, namely, any type internal in a SU-rank 1 type is
s-essentially 1-based (a strong version of the notion "essentially 1-based" from [S1]) by means of
the forking-topology. We do this by generalizing the dichotomy theorem from [S1] to any hypersimple
theory (rather than a countable one) equipped with a projection-closed family of topologies, while
its conclusion is strengthened to get that any type internal in a SU-rank 1 type is s-essentially
1-based (in [S1] we got only "essentially 1-based" in the conclusion), provided that no unbounded
open supersimple is interpreted. This will ensure the existence of many stable formulas that
witness forking. In [S1] we dealt with the remaining case by the development of a model theoretic
Baire category theorem in which we analyze more complicated "forking sets" that are related to the
forking topology. This theorem made an essential use of the existence of many stable formulas and
the assumption that the language is countable.

We assume basic knowledge of simple theories; a good textbook on simple theories is [W]. Throughout
this paper we work in a $\kappa$-saturated and $\kappa$-strongly saturated model $\CC$, for some
large $\kappa$, of a complete first order theory $T$.

\section{The dichotomy}
In this section we assume $T=T^{eq}$ is a hypersimple theory and we work in $\CC=\CC^{eq}$. First
recall the definition of a projection-closed family of topologies.

\begin{definition}\em
A family $$\Upsilon=\{\Upsilon_{x,A} \vert\ x \mbox{ is a finite sequence of variables and }
A\subset \CC \mbox{ is small}\}$$ is said to be \em a projection-closed family of topologies \em if
each $\Upsilon_{x,A}$ is a topology on $S_x(A)$ that refines the Stone-topology on $S_x(A)$, this
family is invariant under automorphisms of $\CC$ and change of variables by variables of the same
sort, the family is closed under product by the full Stone spaces $S_y(A)$ (where $y$ is a disjoint
tuple of variables) and closed by projections, namely whenever $\UU(x,y)\in \Upsilon_{xy,A}$,
$\exists y\UU(x,y)\in\Upsilon_{x,A}$.
\end{definition}

From now on $\Upsilon$ denotes a projection-closed family of topologies.

\begin{definition}\label {def ess-1-based}\em
1) A type $p\in S(A)$ is said to be \em s-essentially 1-based over $A_0\subseteq A$ by means of
$\Upsilon$ \em if for every finite tuple $\bar c$ from $p$ and for every $\Upsilon$-open set $\UU$
over $A\bar c$, with the property that $a$ is independent from $A$ over $A_0$ for every $a\in \UU$,
the set $\{a\in \UU \vert\ Cb(a/A\bar c)\not\in bdd(aA_0)\}$ is nowhere dense in the Stone-topology
of $\UU$. We say $p\in S(A)$ is \em s-essentially 1-based by means of $\Upsilon$ \em if $p$ is
s-essentially 1-based over $A$ by means of $\Upsilon$.\\ 2) Let $V$ be an $A_0$-invariant set and
let $p\in S(A_0)$. We say that $p$ is \em analyzable in $V$ by s-essentially 1-based types by means
of $\Upsilon$ \em if there exists $a\models p$ and there exists a sequence $(a_i\vert\
i\leq\alpha)\subseteq dcl(A_0a)$ with $a_\alpha=a$ such that $tp(a_i/A_0\cup\{a_j\vert j<i\})$ is
$V$-internal and s-essentially 1-based over $A_0$ by means of $\Upsilon$ for all $i\leq\alpha$.
\end{definition}

In [S1] we said that \em $p\in S(A)$ is essentially 1-based with respect to $\Upsilon$, \em if 1)
in Definition \ref{def ess-1-based} holds with the additional requirement that $\UU$ is
type-definable. Before stating the main theorem, recall that for an $A$-invariant set $\UU$ and a
type $p$ over $A$, we say that \em $\UU$ is almost $p$-internal (over $A$) \em if $tp(a/A)$ is
almost $p$-internal for every $a\in\UU$. Also, $\UU$ is said to be \em unbounded \em if it contains
the solution set of some non-algebraic type (equivalently, its cardinality is $\geq\kappa$). We can
now phrase the dichotomy.

\begin{theorem}\label{dichotomy thm}
Let $T$ be any hypersimple theory. Let $\Upsilon$ be a projection-closed family of topologies. Let
$p_0$ be a partial type over $\emptyset$ of $SU$-rank 1. Then, either there exists an unbounded
$\Upsilon$-open set (over some small set $A$) that is almost $p_0$-internal (and in particular has
finite $SU$-rank ), or every complete type $p\in S(A)$ that is internal in $p_0$ is s-essentially
1-based over $\emptyset$ by means of $\Upsilon$. In particular, either there exists an unbounded
$\Upsilon$-open set that is almost $p_0$-internal, or whenever $p\in S(A)$ and every non-algebraic
extension of $p$ is non-foreign to $p_0$, $p$ is analyzable in $p_0$ by s-essentially 1-based types
by means of $\Upsilon$.
\end{theorem}

Before proving the dichotomy, note the following easy generalization of [S1, Proposition 4.4]
(recall the domination notation: $b\unrhd_a c$ iff for any $d$ if $d$ independent from $b$ over $a$
then $d$ is independent from $c$ over $a$.)

\begin{proposition}\label {open Cb}
Let $q(x,y)\in S(\emptyset)$ and let $\chi(x,y,z)$ be an $\emptyset$-invariant set such that for
all $(c,b,a)\models \chi(x,y,z)$ we have $b\unrhd_a bc$. Then the set $$\UU=\{(e,c,b,a) \vert\ e\in
acl(Cb(cb/a))\}$$ is relatively Stone-open inside the set
$$F=\{(e,c,b,a)\vert\ \nonforkempty{b}{a}, \models\chi(c,b,a), tp(cb)=q\}.$$ (where
$e$ is taken from a fixed sort too).
\end{proposition}

The proof of Proposition \ref{open Cb} is the same as in [S1], we write it for completeness. Let us
recall the basic notion and fact that are needed for the proof. Recall that a set $\UU$ is said to
be a basic $\tau^f_*$-open set over $C$ if there exists $\psi(x,y,C)\in L(C)$ such that $\UU=\{a
\vert \ \psi(x,aC) \mbox{\ forks\ over } a\}$.

\begin{fact}\label {tau_star}[S1, Lemma 4.3]
Let $C$ be any set and let $\WW=\{(e,a) \vert\ e\in acl(Cb(C/a))\}$ (where $e,a$ are taken from
fixed sorts). Then $\WW$ is a $\tau^f_*$-open set over $C$.
\end{fact}

\proof Note that since $q\in S(\emptyset)$, it is enough to show that for any fixed $c^*b^*\models
q$ the set $\UU^*=\{(e,a) \vert\ e\in acl(Cb(c^*b^*/a))\}$ is relatively Stone-open inside
$$F^*=\{(e,a)\vert\ \nonforkempty{b^*}{a}, \models\chi(c^*,b^*,a)\}.$$ Now, by Fact \ref{tau_star},
we know $\UU^*$ is a $\tau_*^f$-open set over $b^*c^*$. Thus, for some $\psi_i(t_i;w,z,c^*b^*)\in
L(c^*b^*)$ ($i\in I$) we have $\UU^*=\bigcup_i \UU^*_{\psi_i}$ where $$\UU^*_{\psi_i}=\{(e,a)\vert\
\psi_i(t_i;e,a,c^*b^*) \mbox{\ forks over } ea\}.$$

\begin{subclaim}\label{subclaim 1}
For every $(e,a)\in F^*$ we have $(e,a)\in \UU^*_{\psi_i}$ iff $$\forall d
(\psi_i(d;e,a,c^*b^*)\rightarrow \forkempty{da}{b^*})\wedge e\in acl(a).$$
\end{subclaim}

\proof Let $(e,a)\in F^*$. Assuming the left hand side we know $e\in acl(Cb(c^*b^*/a))$, hence
$e\in acl(a)$. Let $d\models\psi_i(z;e,a,c^*b^*)$. If $\nonforkempty{da}{b^*}$, then
$\nonfork{d}{b^*}{a}$. Since $(e,a)\in F^*$, $b^*\unrhd_a b^*c^*$ implies
$\nonfork{d}{b^*c^*}{ea}$, contradicting $(e,a)\in \UU^*_{\psi_i}$. Assume now the right hand side.
By a way of contradiction assume there exists $d\models\psi_i(t_i;e,a,c^*b^*)$ such that
$\nonfork{d}{b^*c^*}{ea}$. Since $e\in acl(a)$, this is equivalent to $\nonfork{d}{b^*c^*}{a}$.
Since
$(e,a)\in F^*$ this is equivalent to $\nonforkempty{da}{b^*}$, contradiction.$\ \ \ \ \Box$\\

\noindent By Subclaim \ref{subclaim 1} we see that each of $\UU^*_{\psi_i}$ and hence $\UU^*$ is
Stone-open relatively inside $F^*$ (since dependence in $b^*$ is a Stone-open condition over $b^*$).$\ \ \ \ \Box$\\

\noindent\textbf{Proof of Theorem \ref{dichotomy thm}} $\Upsilon$ will be fixed and we'll freely
omit the phrase "by means of $\Upsilon$". To see the "In particular" part, work over $A$ and assume
that every $p'\in S(A')$, with $A'\supseteq A$, that is internal in $p_0$, is s-essentially 1-based
over $A$. Moreover, assume $p\in S(A)$ is non-algebraic and every non-algebraic extension of $p$ is
non-foreign to $p_0$. Then, for $a\models p$ there exists $a'\in dcl(Aa)\backslash acl(A)$ such
that $tp(a'/A)$ is $p_0$-internal and thus s-essentially 1-based over $A$ by our assumption. Thus,
by repeating this process we get that $p$ is analyzable in $p_0$ by s-essentially 1-based types. We
now prove the main part. Assume there exists $p\in S(A)$ that is internal in $p_0$, and $p$ is not
s-essentially 1-based over $\emptyset$. By the definition, there exist a finite tuple $d$ of
realizations of $p$ and $b$ that is independent from $d$ over $A$, and a finite tuple $\bar
c\subseteq p_0$ such that $d\in dcl(Ab\bar c)$, and there exists a $\Upsilon$-open set $\UU$ over
$Ad$ such that $a$ is independent from $A$ for all $a\in \UU$ and $\{a\in \UU \vert
Cb(a/Ad)\not\subseteq acl(a)\}$ is not nowhere dense in the Stone-topology of $\UU$. So, since
$\Upsilon$ refines the Stone-topology, by intersecting $\UU$ with a definable set, we may assume
that $\{a\in \UU \vert Cb(a/Ad)\not\subseteq acl(a)\}$ is dense in the Stone-topology of $\UU$.
\noindent Now, for each (finite) subsequence $\bar c_0$ of $\bar c$, let $$F_{\bar c_0}=\{ a\in \UU
\vert\ \exists b',\bar c'_0,\bar c'_1\ \mbox{s.t.}\ tp(b'\bar c'_0\bar c'_1/Ad)=tp(b\bar c_0(\bar
c\backslash \bar c_0)/Ad)\ \mbox{and} \nonforkempty{a}{Ab'\bar c'_0}\}.$$ Note that since $d$ is
independent from $b$ over $A$, any $a\in\UU$ is independent from $Ab'$ whenever
$tp(b'/Ad)=tp(b/Ad)$ and $\nonfork{a}{b'}{Ad}$. Thus $F_{\langle\rangle}=\UU$. Let $\bar c^*_0$ be
a maximal subsequence (with respect to inclusion) of $\bar c$ such that $F_{\bar c^*_0}$ has
non-empty Stone-interior in $\UU$ over $Ad$ (note that $F_{\bar c}$ has no Stone-interior
relatively in $\UU$). Let $\UU^*=\bigcap_{\bar c^*_0\subset\bar c'\subseteq\bar c} \UU\backslash
F_{\bar c'}$. Note that each $F_{\bar c'}$ is Stone closed relatively in $\UU$. Thus $\UU^*$ is
Stone-dense and open in $\UU$ and therefore there exists a
non-empty relatively Stone-open in $\UU$ set $W^*\subseteq F_{\bar c_0^*}\cap \UU^*$.

\begin{subclaim}\label{subclaim0_main}
$W^*$ is a non-empty $\Upsilon$-open set over $Ad$ such that $\{a\in W^* \vert\
Cb(a/Ad)\not\subseteq acl(a)\}$ is dense in the Stone-topology of $W^*$ and for every $a\in W^*$ we
have: there exists $b'\bar c'_0\bar c'_1\models tp(b\bar c^*_0(\bar c\backslash \bar c^*_0)/Ad)$
such that $a$ is independent from $Ab'\bar c'_0$ over $\emptyset$ and moreover, for every $b'\bar
c'_0\bar c'_1\models tp(b\bar c^*_0(\bar c\backslash \bar c^*_0)/Ad)$ such that $a$ is independent
from $Ab'\bar c'_0$ we necessarily have $\bar c'_1\in acl(aAb'\bar c'_0)$.
\end{subclaim}

\proof As $p_0$ has $SU$-rank 1, this is a conclusion of our construction.$\ \ \ \ \Box$\\

\noindent Let us now define a set $V$ over $Ad$ by\\  $$V=\{(e',b',\bar c'_0,\bar c'_1,a') \vert \
\mbox{if}\ tp(b'\bar c'_0\bar c'_1/Ad)=tp(b\bar c^*_0(\bar c\backslash \bar c_0^*)/Ad)\ \mbox{and}
\nonforkempty{a'}{Ab'\bar c'_0}$$  $$\mbox{then}\ e'\in acl(Cb(Ab'\bar c'_0\bar c'_1/a'))\}.$$

\noindent Let $V^*=\{e' \vert \exists a'\in W^*\ \forall b',\bar c'_0,\bar c'_1\ V(e',b',\bar
c'_0,\bar c'_1,a')\}.$

\begin{subclaim}\label{subclaim1_main}
$V^*$ is a $\Upsilon$-open set over $Ad$.
\end{subclaim}

\proof By Proposition \ref{open Cb} and Subclaim \ref{subclaim0_main}, there exists a Stone-open set $V'$
over $Ad$ such that for all $a'\in W^*$ and for all $e',b',\bar c'_0,\bar c'_1$ we have
$V'(e',b',\bar c'_0,\bar c'_1,a')$ if and and only if $V(e',b',\bar c'_0,\bar c'_1,a')$. Thus, we
may replace $V$ by $V'$ in the definition of $V^*$. As Stone-open sets are closed under the
$\forall$ quantifier, the $\Upsilon$ topology refines the Stone-topology and closed under product
by a full Stone-space and closed under projections, we conclude that $V^*$ is a $\Upsilon$-open
set.$\ \ \ \ \Box$

\begin{subclaim}\label{subclaim2_main}
For appropriate sort for $e'$, the set $V^*$ is unbounded and is almost $p_0$-internal (over $Ad$)
and thus has finite $SU$-rank over $Ad$.
\end{subclaim}

\proof First, note the following general observation.

\begin{fact}\label{dcl_cb remark}\em
Assume $d\in dcl(c)$. Then $Cb(d/a)\in dcl(Cb(c/a))$ for all $a$.
\end{fact}

\noindent Let $a^*\in W^*$ be such that $Cb(a^*/Ad)\not\subseteq acl(a^*)$. Then
$Cb(Ad/a^*)\not\subseteq acl(Ad)$. By Fact \ref{dcl_cb remark}, there exists $e^*\not\in acl(Ad)$
such that $e^*\in acl(Cb(Ab'\bar c'_0\bar c'_1/a^*))$ for all $b'\bar c'_0\bar c'_1\models tp(b\bar
c^*_0(\bar c\backslash \bar c^*_0)/Ad)$. In particular, $e^*\in V^*$. Thus, if we fix the sort for
$e'$ in the definition of $V^*$ to be the sort of $e^*$, then $V^*$ is unbounded. Now, let $e'\in
V^*$. Then for some $a'\in W^*$, $\models V(e',\bar c'_0,\bar c'_1,b',a')$ for all $b',\bar
c'_0,\bar c'_1$. By Subclaim \ref{subclaim0_main}, there exists $b'\bar c'_0\bar c'_1\models
tp(b\bar c^*_0(\bar c\backslash \bar c^*_0)/Ad)$ such that $a'\ \mbox{is\ independent\ from}\
Ab'\bar c'_0\ \mbox{over}\ \emptyset$. Thus, by the definition of $V^*$ and $V$, $e'\in
acl(Cb(Ab'\bar c'_0\bar c'_1/a'))$. Since $Ab'$ is independent from $a'$ over $\emptyset$, $tp(e')$
is almost-$p_0$-internal (as $Cb(Ab'\bar c'_0\bar c'_1/a')$ is in the definable closure of any
Morley sequence of $Lstp(Ab'\bar c'_0\bar c'_1/a')$ ), and in particular $tp(e'/Ad)$ is almost
$p_0$-internal (note that, in general, whenever $q=tp(a/A)$ is internal in an $\emptyset$-invariant
set $\RR$
then any extension of $q$ is almost $\RR$-internal) and therefore $tp(e'/Ad)$ has finite $SU$-rank. $\ \ \ \ \Box$\\

\noindent Thus $V^*$ is the required set. $\ \ \ \ \Box$\\

We now draw some consequences of the above dichotomy for countable languages.

\begin{theorem}\label{cor1}
Let $T$ be any countable hypersimple theory. Let $\Upsilon$ be a projection-closed family of
topologies such that $\{ a\in \CC^x \vert a\not\in acl(A)\}\in \Upsilon_{x,A}$ for all $x$ and set
$A$ . Let $p_0$ be a partial type over $\emptyset$ of $SU$-rank 1. Then, either there exists an
unbounded type-definable $\Upsilon$-open set over some small set that is almost $p_0$-internal and
has \textbf{bounded} finite $SU$-rank, or every complete type $p\in S(A)$ that is internal in $p_0$
is essentially 1-based over $\emptyset$ by means of $\Upsilon$. In particular, either there exists
an unbounded $\Upsilon$-open set that is almost $p_0$-internal and has \textbf{bounded} finite
$SU$-rank, or whenever $p\in S(A)$, where $A$ is countable, and every non-algebraic extension of
$p$ is non-foreign to $p_0$, $p$ is analyzable in $p_0$ by essentially 1-based types by means of
$\Upsilon$.
\end{theorem}

\proof We go back to the proof of Theorem \ref{dichotomy thm} (the main part); we start with $p\in
S(A)$ that is $p_0$-internal and not essentially 1-based over $\emptyset$ and apply the same proof
(but note that in the proof of Theorem \ref{dichotomy thm} we assumed $p$ is not s-essentially
1-based). So, now $\UU$ is assumed to be a \textbf{type-definable} $\Upsilon$-open set over $Ad$.

\begin{subclaim}
We may assume $W^*$ is type-definable and $\Upsilon$-open over $Ad$ and there exists
$V^{**}\subseteq V^*$ that is unbounded, type-definable and $\Upsilon$-open over $Ad$.
\end{subclaim}

\proof In the proof of Theorem \ref{dichotomy thm} the set $W^*$ is chosen to be a non-empty
intersection of $\UU$ with a Stone-open set over $Ad$, so we could instead take it to be a
non-empty intersection of $\UU$ with a definable subset of this Stone-open set (and still
$W^*\subseteq F_{\bar c_0^*}\cap \UU^*$). Since $\UU$ is $\Upsilon$-open and type-definable, $W^*$
is type-definable and $\Upsilon$-open over $Ad$. Now, by the definition of $V^*$ and the proof of
Subclaim \ref{subclaim1_main} there exist a Stone open set $V_0$ over $Ad$ such that $V^*=\{ e'
\vert \exists a'\in W^*\ (V_0(e',a'))\}$. From this we easily get the required set $V^{**}$ (by
replacing $V_0$ by a definable set and using the fact that $W^*$ is type-definable and that $\Upsilon$ is a
projection-closed family of topologies). $\ \ \ \ \Box$\\

\noindent By the proof of Subclaim \ref{subclaim2_main} we know that for all $e'\in V^{**}$ we have
$e'\in acl(Cb(Ab'\bar c'_0\bar c'_1/a'))$ for some $a'\in W^*$ and some $b',\bar c'_0,\bar c'_1$
such that $a'$ is independent from $Ab'\bar c'_0\ \mbox{over}\ \emptyset$ and $b'\bar c'_0\bar
c'_1\models tp(b\bar c^*_0(\bar c\backslash \bar c^*_0)/Ad)$. Let $q=tp(Ab\bar c^*_0)$. For every
$\chi=\chi(x,y_0,...,y_n,\bar z)\in L$ (for some $n<\omega$) such that $\forall y_0 y_1 ...y_n \bar
z\ \exists^{<\infty}x\ \chi(x,y_0,y_1,...y_n,\bar z)$, and $m<\omega$ let

$$F_{\chi,m}=\{e\in V^{**} \vert\ \models\chi(e,C_0,C_1,..C_n,\bar c)\ \mbox{for\ some\ } \bar c\in
p_0^m\ \mbox{and\ some}\ \emptyset-\mbox{independent\ sequence}$$ $\ \ \ \ \ \ \ \ \ \ \ \ \ \ \ \
\ \ \ \ \ \ \ \ \ (C_i\vert i\leq n)\ \mbox{of\ realization of}\ q\ \mbox{with}
\nonforkempty{e}{(C_i\vert i\leq n)} \}.$\\

\noindent By the aforementioned, we get that $V^{**}\subseteq \bigcup_{m,\chi} F_{m,\chi}$ (the
union is over each $m,\chi$ as above). By the Baire category theorem applied to the Stone-topology
of the Stone-closed set $V^{**}\backslash acl(Ad)$, there exists $\theta\in L(Ad)$ such that
$$\tilde V\equiv \theta^\CC\cap(V^{**}\backslash acl(Ad))\neq\emptyset\ \mbox{and\ } \tilde
V\subseteq F_{m^*,\chi^*}$$ for some $m^*,\chi^*$ as above. Clearly, $\tilde V$ is unbounded,
type-definable and $\Upsilon$-open (by the assumptions on $\Upsilon$). Now, for every $a\in \tilde
V$, $SU(a/Ad)\leq m^*$ and $tp(a/Ad)$ is almost $p_0$-internal (as $tp(a)$ is almost
$p_0$-internal, and $SU(a)\leq m^*$ by the definition of $F_{m^*,\chi^*}$). This completes the
proof of the first part of the theorem. The rest follow easily by repeated applications of the
first part (when working over $A$).$\ \ \ \ \Box$

$\\$

Recall that $T$ is PCFT if its forking-topologies is a projection-closed family of topologies, that
is, whenever $\UU(x,y)$ is a $\tau^f$-open set over a small set $A$, $\exists y\UU(x,y)$ is a
$\tau^f$-open set over $A$. Applying Theorem \ref{cor1} for the special case of the
forking-topologies we conclude the following.

\begin{corollary}
Let $T$ be any countable hypersimple theory with PCFT. Let $p_0$ be a partial type over $\emptyset$
of $SU$-rank 1. Then, either there exists a weakly-minimal formula that is almost $p_0$-internal,
or every complete type $p\in S(A)$ that is internal in $p_0$ is essentially 1-based over
$\emptyset$ by means of $\tau^f$. In particular, either there exists a weakly-minimal formula that
is almost $p_0$-internal, or whenever $p\in S(A)$, where $A$ is countable, and every non-algebraic
extension of $p$ is non-foreign to $p_0$, $p$ is analyzable in $p_0$ by essentially 1-based types
by means of $\tau^f$.

\end{corollary}

\proof Our assumptions are clearly a special case of the assumptions of Theorem \ref{cor1}, thus
we only need to prove the first part. By the conclusion of Theorem \ref{cor1}, we may assume that there
exists a $\tau^f$-open set $\UU$ of bounded finite $SU$-rank over some small set $A$ that is almost
$p_0$-internal. Recall now [S0, Proposition 2.13]:

\begin{fact}\label{tau bounded SU}
Let $\UU$ be an unbounded $\tau^f$-open set over some set $A$. Assume $\UU$ has bounded finite
$SU$-rank. Then there exists a set $B\supseteq A$ and $\theta(x)\in L(B)$ of $SU$-rank 1 such that
$\theta^\CC\subseteq \UU\cup acl(B)$.
\end{fact}

By Fact \ref{tau bounded SU}, there exists exists a weakly-minimal $\theta(x,b)\in L(B)$ for some
small set $B\supseteq A$, such that $\theta^\CC\subseteq \UU\cup acl(B)$. Now, $tp(a/B)$ is almost
$p_0$-internal for every $a\in \theta^\CC$, and so $tp(a/b)$ ($b$ is the parameter of $\theta(x,b)$
) is almost $p_0$-internal over $b$ for every $a\in \theta^\CC$ (by taking non-forking
extensions).$\ \ \ \ \Box$\\

We now state the main conclusion for uncountable hypersimple unidimensional theories.

\begin{definition}
We say that $T$ is s-essentially 1-based if for every $SU$-rank 1 partial type $p_0$ over some $A$,
every $p\in S(A)$ that is internal in $p_0$ is s-essentially 1-based by means of $\tau^f$.
\end{definition}

\begin{corollary}\em
Let $T$ be a hypersimple unidimensional theory that is not s-essentially 1-based. Then $T$ is
supersimple.
\end{corollary}

\proof First, recall the following fact [S1, Corollary 3.15] (an $A$-invariant set $\UU$ is called
\em supersimple \em if $SU(a/A)<\infty$ for every $a\in\UU$).

\begin{fact}\label {fact1}
Let $T$ be a hypersimple unidimensional theory and work in $\CC=\CC^{eq}$. Let $p\in S(A)$ and let
$\UU$ be an unbounded $\tau^f$-open set over $A$. Then $p$ is analyzable in $\UU$ in finitely many
steps. In particular, for such $T$ the existence of an unbounded supersimple $\tau^f$-open set over
some small set $A$ implies $T$ is supersimple.
\end{fact}

\noindent Now, assume $T$ is a hypersimple unidimensional theory that is not s-essentially 1-based.
By Theorem \ref{dichotomy thm}, there exists an unbounded $\tau^f$-open set of finite $SU$-rank
over some small set. By Fact \ref{fact1}, every complete type has finite $SU$-rank.$\ \ \ \ \Box$

\end{document}